\documentstyle[12pt]{article}

\newcommand{\sect}[1]{\section{#1}}

\newtheorem{theorem}{Theorem}[section]
\newtheorem{proposition}[theorem]{Proposition}
\newtheorem{definition}[theorem]{Definition}

\newtheorem{lemma}[theorem]{Lemma}
\newtheorem{remark}[theorem]{Remark}

\def\lan{\langle}
\def\ran{\rangle}

\def\pair#1#2{\lan #1,#2\ran}

\def\fidi{\hskip5pt \vrule height4pt width4pt depth0pt \par}

\def\1{\'{\i}}

\def\De{\Delta}

\def\CC{{\bf C}}

\def\ZZ{{\bf Z}}

\def\id{{\rm id}}

\def\half{{\scriptstyle{1\over2}}}

\def\tens{\otimes}
\def\fraz#1#2{{\strut\displaystyle #1\over\displaystyle #2}}

\def\indl{{\rm ind}_K^G(\rho_L)}
\def\indr{{\rm ind}_K^G(\rho_R)}

\def\fq{{\cal F}_q}
\def\fqg{{\cal F}_q(G)}
\def\fqk{{\cal F}_q(K)}
\def\ffq#1{{\cal F}_q(#1)}

\def\uq{{\cal U}_q}
\def\rin{\rho_{\cal I}}

\def\urin{\widetilde\rho_{\cal I}}
\def\urinn#1{\widetilde\rho_{{\cal I}_{#1}}}
\def\ulin{\widetilde\lambda_{\cal I}}
\def\ulinn#1{\widetilde\lambda_{{\cal I}_{#1}}}
\def\rrmu{{\widetilde\rho_{m,u}}}
\def\th{{\widehat{t}}}
\def\muh{{\widehat{\mu}}}
\def\xh{{\widehat{x}}}

\def\fmu{\varphi_{m,u}}
\def\hirr{{H}_{m,u}^{\,{\rm irr}}}
\def\hoirr{{H}_{0}^{\,{\rm irr}}}
\def\R{{\bf R}}

\begin{document}

\title{Unitarity of induced representations\\
from coisotropic quantum subgroups}
\author{F.Bonechi$\,{}^{1,2}$, N.Ciccoli$\,{}^{3}$, R.Giachetti$\,{}^{2,1}$, 
E.Sorace$\,{}^{1,2}$, M.Tarlini$\,{}^{1,2}$}

\maketitle
\centerline{{\small  ${ }^1$ INFN Sezione di Firenze}}
\centerline{{\small ${ }^2$ Dipartimento di Fisica, Universit\`a di Firenze, Italy. }}
\centerline{{\small ${ }^3$ Dipartimento di Matematica, Universit\`a di Perugia, Italy. }}

\begin{abstract}
{We study unitarity of the induced representations from coisotropic quantum subgroups
which were introduced in \cite{N2}.
We define a real structure on coisotropic subgroups which determines an involution
on the homogeneous space. We give general invariance properties for functionals
on the homogeneous space which are sufficient to build a unitary representation starting
from the induced one. We present the case of the one-dimensional quantum Galilei group,
where we have to use in all generality our definition of quasi-invariant functional.
}
\end{abstract}

\medskip
{\bf Math.Subj.Classification}: 16W30, 17B37 81R50


\thispagestyle{empty}
\sect{Introduction}
\bigskip
Induced representations are the fundamental tool in the representation
theory of Lie groups. Starting with
a unitary representation of a generic closed subgroup $H$ of
a given Lie group $G$ there is a standard way to induce 
a unitary representation of the whole group, realizing it on
the space of $L^2$ sections of a homogeneous vector bundle on $G/H$.
Whenever the homogeneous space is compact, the unitarity 
of the induced representation is due to the existence of an 
invariant measure on it.

For a non compact homogeneous space such an invariant measure 
does not necessarily exists. Nevertheless one can fix a 
quasi invariant measure $\mu$, {\it i.e.} a measure equivalent 
to its translates. A quasi invariant measure
always exists and is unique only up to equivalence.
Once such a measure has been fixed, it is possible to modify 
accordingly the definition of induced representation, 
using a weight function, to obtain a
unitary representation. It is finally possible to
show that equivalent 
measures give rise to equivalent unitary representations.

For quantum groups induction of representations was up to now
analyzed either in the algebraic \cite{PW} or in the compact case \cite{GZ}.
In both cases, the authors start from corepresentations of quantum subgroups.
However it was realized from the beginning that quantum subgroups
are very rare. The more general notion of quantum coisotropic
subgroup is sufficient to obtain all embeddable quantum homogeneous 
spaces. The semiclassical limit of such subgroups has a natural 
interpretation in terms of Poisson
structures: indeed one can observe that
the subgroup obtained by conjugation of a Poisson Lie subgroup is no more
Poisson Lie but only Poisson coisotropic and the conjugated subgroups 
determine inequivalent Poisson structures on the corresponding
homogeneous space, see \cite{N1}. After quantization, this situation gives 
rise to a family of inequivalent quantum homogeneous spaces, which are not
obtained by a quantum subgroup: it is therefore rather natural to induce 
representations starting from quantum coisotropic subgroups. This
program was started in \cite{N2}, where the basic algebraic properties
were proved. 

The purpose of the present work is to study unitarity 
in the same general setting.
We propose a definition of quasi invariant functional on a quantum homogeneous
space which is quite general. We then fix a quasi invariant functional
and we prove a general procedure to define unitary representations
from the induced ones. 

The exposition given in this paper is mainly based on the coalgebra 
properties: strictly speaking such an approach is completely
suited only for the compact case.
The lack of a general definition of non compact quantum group
and quantum homogeneous space does not allow to develop
a meaningful general theory. Nevertheless many concrete
examples exist, in which unitary representations on Hilbert spaces are
built starting from $\uq$-module $*$-algebras, like, for example, the
$E_q(2)$ case in \cite{Koe}, $SU_q(1,1)$ in \cite{No} and the 
quantum Galilei group in \cite{F1}.
We will explicitly study the induced representations of the quantum
Galilei Group showing how our results extend to $\uq$-module $*$-algebras.
In \cite{F1} physical motivations lead to study its unitary 
representations with respect to a non standard involution.
When the involution is chosen to be standard, as we will do
here, induced representations must be studied using a functional 
which is quasi invariant according to our definition.

The paper is structured as follows. In section 2 we introduce coisotropic
quantum subgroups, focusing our attention on reality conditions, which
were only briefly mentioned in \cite{N1}. 
This is necessary to give an involution
on the homogeneous space and to define unitary corepresentations for
the subgroup.

In section 3 we define quasi invariant functionals on homogeneous spaces.
We will call {\it essentially invariant} those functionals already defined 
in \cite{Ko} and we show that an essentially invariant functional
is always equivalent to an invariant one. 

In section 4 we define a sesquilinear form on the representation 
space and we show that there is one definite 
choice according to the fact
that the coisotropic subgroup is left or right (see Lemma \ref{lem_sesq}). 
We then construct the unitary
representation with respect to this sesquilinear form and we finally 
show that equivalent quasi invariant functionals give rise to 
equivalent representations.

In section 5 the example of quantum Galilei group $\uq(G(1))$ is worked out.
Dealing with a non compact quantum homogeneous space we are forced to describe 
it as a $\uq$-module $*$-algebra. We construct its unitary representations by
means of a functional which is quasi invariant but not essentially invariant.

\bigskip
\bigskip
\sect{Real Coisotropic Subgroups}
\bigskip
In this section we will give the definition and main properties of real
coisotropic
quantum subgroups and their associated homogeneous spaces. 
With respect to \cite{N1},
to which we refer for a more detailed discussion, the emphasis here is on 
the $*$
structure on the subgroup and on the corresponding homogeneous space.

\smallskip
\begin{definition}
\label{def_cos}
Given a real quantum group $(\fqg,*)$ we will call {\rm real coisotropic quantum right
(left) subgroup} $(\fqk, \tau_K)$ a coalgebra, right (left) $\fqg$-module $\fqk$ such that:
\begin{itemize}
\item[\rm{($i$)}] there exists a surjective linear map $\pi: \fqg\rightarrow\fqk$, which is a
morphism of $\fqg$-modules (where $\fqg$ is considered as a module on itself via multiplication)
and of coalgebras;
\item[\rm{($ii$)}] there exists an antilinear map $\tau_K: \fqk\rightarrow\fqk$ such that 
$\tau_K\circ\pi=\pi\circ\tau$, where $\tau = *\circ S$.
\end{itemize}
A $*$-Hopf algebra $\fqk$ is said to be a real quantum subgroup if there exists a 
$*$-Hopf algebra epimorphism $\pi:\fqg\rightarrow\fqk$.
\end{definition}

\medskip
\begin{remark} {\rm A coisotropic quantum subgroup is not a $*$-coalgebra 
but it has only $\tau_K$ defined on it. It is easy to verify that if a coisotropic 
quantum subgroup is also a
$*$-coalgebra and $\pi\circ *=*\circ\pi$, then it is possible to complete the structure
so to have a quantum subgroup.}
\end{remark}

\smallskip
In the following we will always use subgroup to mean real subgroup.
Coisotropic quantum subgroups are easily characterized by the following proposition.
\medskip
\begin{proposition}
\label{pro_coiso}
There exists a bijective correspondence between coisotropic quantum right (left) subgroups 
and $\tau$-invariant two-sided coideals, right (left) ideals in $\fqg$. \fidi
\end{proposition}

\medskip
Coisotropic quantum subgroups give in a canonical way embeddable quantum homogeneous spaces.
\medskip
\begin{definition}
\label{def_emb_hom}
A $*$-algebra $B$ is said to be an embeddable quantum left (right) $\fqg$-homogeneous space 
if there exists a coaction $\delta: B\rightarrow B\otimes \fqg$, 
($\delta: B\rightarrow \fqg\otimes B$) and an injective
morphism of $*$-algebras $i:B\rightarrow\fqg$ such that 
$\Delta\circ i = (i\otimes \id)\circ\delta$ ($\Delta\circ i = (\id\otimes i)\circ\delta$).
\end{definition}
\medskip
In the following we will identify left (right) embeddable quantum homogeneous spaces with 
$*$-subalgebras and right (left) coideals of $\fqg$. 
\medskip

\begin{proposition}
\label{pro_emb_hom} 
Let $(\fqk,\tau_K)$ be a left coisotropic quantum subgroup of 
$(\fqg$, $*)$ then
$$
B_\pi = \{a\in\fqg \,|\, (\pi\otimes\id)\Delta a = \pi(1)\otimes a\}
$$
is a left embeddable quantum homogeneous space. 

If $(\fqk,\tau_K)$ is a right coisotropic quantum 
subgroup of ($\fqg$,$*$) then
$$
B^\pi = \{a\in\fqg \,|\, (\id\otimes\pi)\Delta a = a\otimes\pi(1)\}
$$
is a right embeddable quantum homogeneous space.
\end{proposition}
\smallskip
{\it Proof}. 
Let $\fqk$ be a left coisotropic quantum subgroup.
As a consequence of the properties of $\pi$ it is easy to see 
that if $a,b\in B_\pi$ then $(\pi\otimes\id)\Delta ab=\pi(1)\otimes ab$ and 
$(\pi\otimes\id\otimes\id)(\Delta\otimes\id)\Delta a = \pi(1)\otimes\Delta(a)$,
which means that $B_\pi$ is a subalgebra and right coideal of $\fqg$. 
We have to show that $B_\pi$ is also $*$-invariant. 

If $a\in B_\pi$, then 
$$
\begin{array}{rcl}
\epsilon(a) \pi(1) &=& \sum_{(a)} \pi(S^{-1}(a_{(2)})a_{(1)}) = 
\sum_{(a)} S^{-1}(a_{(2)})\pi(a_{(1)}) \cr
&=& \sum_{(a)} \epsilon(a_{(1)})S^{-1}(a_{(2)}) \pi(1) = \pi(S^{-1}(a)) \,,\cr
\end{array}
$$
where we used the property of $\pi$ with respect to the left module structure of $\fqk$ and
the structure of right coideal of $B_\pi$. 
Applying $\tau_K$ and using commutativity with $\pi$, we conclude that
$$\pi(a^*) = \epsilon(a^*) \pi(1) \;.$$
Because $B_\pi^*$ is a subalgebra and right coideal of $\fqg$ then, using Proposition
$4.7$ of \cite{N1} we have that $(\pi\otimes\id)\Delta a^*=\pi(1)\otimes a^*$ and so
$B_\pi^* =B_\pi$. 

The proof for right coisotropic quantum subgroups is very similar. \fidi

\medskip
\begin{remark}{\rm
A left coisotropic quantum subgroup cannot give rise in general to a right homogeneous 
space (referring to the proof of Proposition \ref{pro_emb_hom}, note that 
$\pi(S(a))\neq \pi(a)$, for $a\in B_\pi$). 
}
\end{remark}

\bigskip
\bigskip
\sect{Quasi-invariant Functionals}
We denote by the same symbol $X.a$ the left regular representation of 
$\uq$ on $\fqg$ and its restriction to the left quantum homogeneous space 
$B_\pi$, {\it i.e.} $X.a=(\id\otimes X)\Delta a\in B_\pi$ if $a\in\ B_\pi$. 
Analogously we denote by 
$a.X$ the restriction of the right regular representation of $\uq$ on $\fqg$ to the
right quantum homogeneous space $B^\pi$, {\it i.e.} 
$a.X=(X\otimes\id)\Delta a\in B^\pi$ if $a\in\ B^\pi$.

\medskip
\begin{definition}
\label{m_equiv}
We say that the linear functional $h:B_\pi\rightarrow{\bf C}$ ({\rm resp.} 
$h:B^\pi\rightarrow{\bf C}$) is real if $h(a^*)={\overline{h(a)}}$; 
furthermore if $h(a^*a)\geq 0$ for
each $a\in B_\pi$ ({\rm resp.} $a\in B^\pi$) $h$ is said to be positive. 
Two such real functionals  $h_1,h_2$ are $\xi$ equivalent if there 
exists an invertible $\xi\in B_\pi$ such that
$$
h_1(a) = h_2(\xi^*a\xi) \;.
$$
\end{definition}

It is easily verified that Definition \ref{m_equiv} defines an equivalence relation.
\medskip
\begin{definition}
\label{def_m_q_inv}
Let $\phi:\uq\rightarrow B_\pi$ be a linear application which satisfies 
the following properties:
\begin{equation}
\label{def_phi}
\phi[XY] = \sum_{(X)} X_{(1)}.\phi[Y] \,\phi[X_{(2)}];\quad\quad 
\phi[1_{\uq}] = 1_{B_\pi} ,
\end{equation}
and $h:B_\pi\rightarrow\CC$ a real functional. We say that
$h$ is {\rm quasi invariant} with weight $\phi$ if 
the following relation is valid for each $X\in\uq$
\begin{equation}
\label{def_h}
h(X.a) = \sum_{(X)} h(\phi[X_{(1)}^*]^*a\phi[S(X_{(2)})])  \;.
\end{equation}
If there exists an invertible $\xi\in B_\pi$ 
such that $\phi[X] = X.\xi \,\xi^{-1}$ then $h$ is said to be 
{\rm essentially invariant}.  
\end{definition}

\medskip
\begin{remark}{\rm
($i$) If $\xi$ is group like and $\tau(\xi)=\xi$ then the essentially invariant
functional with weight $\phi[X]=\langle X, \xi \rangle 1_{B_\pi}$ is invariant, 
{\it i.e.} $h(X.a)=\epsilon(X) h(a)$.

($ii$) Definition \ref{def_m_q_inv} is well posed: indeed it preserves reality of $h$
and is compatible with the regular representation.

($iii$) An essentially invariant functional with a group-like $\xi$ is quasi invariant
according to \cite{Ko}.}
\end{remark}

\medskip

\begin{lemma}
\label{lem_q_inv_m}
A real functional $h:B_\pi\rightarrow \CC$ is quasi-invariant if and only if 
there exists a linear application $\phi:\uq\rightarrow B_\pi$ which satisfies 
{\rm(\ref{def_phi})} and such that, for each $X\in\uq$, we have
\begin{equation}
\label{prop_h}
\sum_{(X)} h(X_{(1)}.a\,\phi[X_{(2)}] ) = h(\phi[X^*]^* a) \;.
\end{equation}
\end{lemma}
\smallskip
{\it Proof}. Let $h$ verify (\ref{prop_h}), then
$$
\begin{array}{rcl}
\sum_{(X)} h(\phi[X_{(1)}^*]^*a\phi[S(X_{(2)})]) &=&
\sum_{(X)} h( X_{(1)}.( a\phi[S(X_{(3)})] ) \phi[X_{(2)}] ) \cr
&=& \sum_{(X)} h(X_{(1)}.a\,X_{(2)}.\phi[S(X_{(4)})] \,\phi[X_{(3)}] ) \cr
&=& h(X.a) \;,
\end{array}
$$
where we applied (\ref{prop_h}) at the first line and (\ref{def_phi}) to obtain the last one.
Let now $h$ be quasi-invariant, then 
$$
\begin{array}{rcl}
\sum_{(X)} h(X_{(1)}.a\phi[X_{(2)}]) &=& 
\sum_{(X)} h(X_{(1)}.(aS(X_{(2)}).\phi[X_{(3)}]) )   \cr
&=& \sum_{(X)} h(\phi[X_{(1)}^*]^*a \,S(X_{(3)}).\phi[X_{(4)}] \,\phi[S(X_{(2)})] )  \cr
&=& h(\phi[X^*]^*a)\;,
\end{array}
$$
where we used trivial properties of the regular representation in the first line, 
relation (\ref{def_h}) in the second one and (\ref{def_phi}) in the last one. \fidi

\bigskip
We give the equivalent definition and lemma (without proof) for right homogeneous spaces.

\medskip 
\begin{definition}
\label{def_m_q_inv2}
Let $\psi:\uq\rightarrow B^\pi$ be a linear application which satisfies the following properties:
\begin{equation}
\label{def_psi}
\psi[XY] = \sum_{(X)} \psi[Y_{(1)}] \psi[X].Y_{(2)} ;\quad\quad 
\psi[1_{\uq}] = 1_{B^\pi} ,
\end{equation}
and $h:B^\pi\rightarrow\CC$ a real functional. We say that
$h$ is {\rm quasi invariant} if the following relation is
valid for each $X\in\uq$
\begin{equation}
\label{def_h2}
h(a.X) = \sum_{(X)} h(\psi[S(X_{(1)})]\,a\,\psi[X_{(2)}^*]^*)  \;.
\end{equation}
The functional $h$ is said to be {\rm essentially invariant} if,
for each $X\in\uq$, $\psi[X]=\xi^{-1}\xi.X$ for some invertible $\xi\in B^\pi$. 
\end{definition}

\medskip
\begin{lemma}
\label{lem_q_inv_m2}
A real functional $h:B^\pi\rightarrow \CC$ is quasi-invariant if and only 
if there exists a linear application $\psi:\uq\rightarrow B^\pi$ which 
satisfies {\rm(\ref{def_psi})} and such that
\begin{equation}
\label{prop_h2}
\sum_{(X)} h(\psi[X_{(1)}]\,a.X_{(2)}) = h(a\psi[X^*]^*) \;.
\end{equation}
\end{lemma}

\medskip
\begin{remark}{\rm
Definitions \ref{def_m_q_inv} and \ref{def_m_q_inv2} don't require the coalgebra 
structure of $B_\pi$ and $B^\pi$. They can be used in the general setting of 
$\uq$-module left and right $*$-algebras.
In the left case they are algebras and left $\uq$-modules such that 
$X.(ab)=\sum_{(X)} X_{(1)}.a \,X_{(2)}.b$, $\,X.1=\epsilon(X)1$ and $(X.a)^*=\tau(X).a^*$;
analogously for the right case.
}
\end{remark}

\medskip
Until the end of this section we will consider only left homogeneous
spaces. Similar statements hold in the right case with obvious
modifications.

In the following Proposition we list some properties of quasi invariant functionals.

\begin{proposition}
\label{m_inv_cor}
Let $B_\pi$ be a left homogeneous space. Then
\begin{itemize}
\item[{\rm($i$)}] if two real functionals $h_i$, $i=1,2$, are $\xi$-equivalent and $h_2$ is
quasi invariant with weight $\phi_2$ then $h_1$ is quasi invariant with weight
\begin{equation}
\label{w_equiv}
\phi_1[X] =\sum_{(X)}X_{(1)}.\xi \phi_2[X_{(2)}]\,\xi^{-1}  \;;
\end{equation}
\item[\rm($ii$)] if $h$ is a quasi invariant real functional with weight $\phi$ 
and $k\in\uq$ such that $\De(k) = k\otimes k$ and $\tau(k)=k$ then the translated
functional $h_k(a)=h(k.a)$ is an equivalent real functional, which is quasi invariant
with weight
$$\phi_k[X] = \phi[XS(k)]\,S(k).\phi[k]$$
and is positive if $h$ is positive;
\item[\rm($iii$)] a real functional is essentially invariant if and only if it is
equivalent to an invariant one.
\end{itemize}
\end{proposition}

\smallskip
{\it Proof}. 
Let us prove item ($i$):
$$
\begin{array}{rcl}
h_1(X.a) &=& h_2(\xi^*\,X.a\,\xi) \\
              &=& \sum_{(X)} h_2(X_{(2)}. 
                  \left(S^{-1}(X_{(1)}).\xi^* \,a\, S(X_{(3)}).\xi\right)) \\
              &=& \sum_{(X)} h_2( (X_{(1)}^*.\xi\phi_2[X_{(2)}^*])^* a 
                  S(X_{(4)}).\xi \,\phi_2[S(X_{(3)})]  ) \\
              &=& \sum_{(X)} h_2((\phi_1[X_{(1)}^*]\xi)^* a \phi_1[S(X_{(2)})]\xi)\\
              &=& \sum_{(X)} h_1(\phi_1[X_{(1)}^*]^* a \phi_1[S(X_{(2)})] )\; .
\end{array} 
$$

Let's prove the item ($ii$). The property for $h_k$ to be a real functional
(and positive if $h$ is positive) directly follows from reality of $k$ 
and from standard duality relations.
>From the property of $h$ of being quasi invariant we have that
$h_k(a)=h(\xi^*a\xi)$ with $\xi=\phi[k^*]$. 
>From the first item
$h_k$ is quasi invariant and its weight is obtained by direct computation.

As a consequence of item ($i$) a functional which is equivalent to an invariant
one is essentially invariant. Let now $h_\xi$ be essentially
invariant with weight $\phi[X]=X.\xi\,\xi^{-1}$. Then always from item ($i$)
$h(a)=h_\xi((\xi^{-1})^*a\xi^{-1})$ is quasi invariant with weight 
$\epsilon(X)$. \fidi

\medskip
\begin{remark}{\rm
($i$) Let $\fqg$ be a compact quantum group according to \cite{DK}. 
On every quantum homogeneous space of a compact quantum group there exists a 
unique invariant real functional $h$ such that $h(1)=1$. 

\smallskip
($ii$) Definitions for weight functions can be rephrased in
a cohomological language. 
Let's define $C^0(B_\pi) = \{\xi\in B_\pi\,|\; \exists \,\xi^{-1} \,\in B_\pi\}$
and $C^k(B_\pi) = \{\phi:\bigotimes^k\uq\rightarrow B_\pi\}$ 
for $k\geq 1$. Then, let $d_0:C^0\rightarrow C^1$, $d_1: C^1\rightarrow C^2$ as
$$
d_0(\xi)[X] = X.\xi\,\xi^{-1} \,,~~~~~~
d_1(\phi)[X\otimes Y] = \phi[XY]-\sum_{(X)} X_{(1)}.\phi[Y] \,\phi[X_{(2)}]  \;.
$$
It is easy to see that $d_1\circ d_0=0$. We can define as usual 
$H^1(B_\pi)={\rm Ker}\,d_1/{\rm Im}\,d_0$. 
If $H^1(B_\pi)=\{0\}$ then every quasi invariant real functional 
is essentially invariant.}
\end{remark}

\bigskip
\bigskip

\sect{Induced Corepresentations}
We summarize in this section general properties of induced corepresentations 
and we start to study their unitarity. We refer to \cite{N2} for a 
detailed study of induced corepresentations from a purely algebraic 
point of view. 

Let ($\fqk,\tau_K$) be a coisotropic quantum left (right) subgroup 
of ($\fqg,*$) and let $B_\pi$ ($B^\pi$) be the associated left (right) quantum 
homogeneous space according to Proposition \ref{pro_emb_hom}. 
Let $\rho_R: V\rightarrow V\otimes\fqk$ be a right corepresentation of 
$\fqk$; if $\{e_i\}$ is an orthonormal basis of $V$ with respect to a
scalar product $\langle,\rangle$ then $\rho_R(e_i) = \sum_j e_j\otimes a_{ji}$. 
Analogously, let $\rho_L: V\rightarrow\fqk\otimes V$ be a left corepresentation 
of $\fqk$, {\it i.e.} $\rho_L(e_i) = \sum_j b_{ij}\otimes e_j$. We say that 
$\rho_R$ ($\rho_L$) is unitary if 
$\tau_K(a_{ij})= a_{ji}$ ($\tau_K(b_{ij})=b_{ji}$).

Let $L=(\pi\otimes\id)\Delta:\fqg\rightarrow\fqk\otimes\fqg$ and 
$R=(\id\otimes\pi)\Delta:\fqg\rightarrow\fqg\otimes\fqk$. We define
$$
\indr = \{F\in V\otimes\fqg \,|\, (\id\otimes L)F = (\rho_R\otimes\id)F\} \;, 
$$
$$
\indl = \{F\in \fqg\otimes V \,|\,(R\otimes \id)F= (\id\otimes\rho_L)F\} \;.
$$

It is easy to prove the following lemma.

\medskip
\begin{lemma}
\label{lem_B_mod}
If $\fqk$ is a left quantum coisotropic subgroup then the linear space $\indr$ is a right
$B_\pi$-module, i.e. if $A=\sum_je_j\otimes A_j\in\indr$, $a\in B_\pi$ then 
$Aa=\sum_j e_j\otimes A_ja\in\indr$. 

If $\fqk$ is a right quantum coisotropic subgroup then the linear space $\indl$ is a left
$B^\pi$-module, i.e. if $A=\sum_jA_j\otimes e_j\in\indl$, $a\in B^\pi$ then 
$aA=\sum_j aA_j\otimes e_j\in\indr$. 
\end{lemma}
\medskip

In \cite{N2} the following proposition is proven.

\medskip
\begin{proposition}
\label{pro_ind_cor}
The restrictions of $(\id\otimes\Delta)$ to $\indr$ and of 
$(\Delta\otimes\id)$ to $\indl$ 
respectively define a right and a left corepresentation of $\fqg$.
\end{proposition}
\medskip
The corepresentations thus defined are called induced 
corepresentations. Induced representations are defined by duality.
To obtain unitary representations we must first find
an invariant sesquilinear form. The following lemma reduces this problem to
that of defining an invariant functional on the corresponding homogeneous space.

\medskip
\begin{lemma}
\label{lem_sesq}
Let $\fqk$ be a coisotropic quantum left subgroup of $\fqg$. Let $\rho_R$ be a unitary
right corepresentation on $V$ and $\langle,\rangle_L$ a 
sesquilinear map from $V\otimes\fqg$ to $\fqg$ given
by $\langle v\otimes a, w\otimes b\rangle_L = \langle v,w\rangle a^*b$. If $A,B \in \indr$
then $\langle A,B \rangle_L \in B_\pi$.

Let $\fqk$ be a coisotropic quantum right subgroup of $\fqg$, $\rho_L$ a unitary left
corepresentation on $V$ and $\langle,\rangle_R$ 
a sesquilinear map from $\fqg\otimes V$ to $\fqg$ given
by $\langle a\otimes v, b\otimes w\rangle_R = \langle v,w\rangle ba^*$. If $A,B \in \indl$
then $\langle A,B \rangle_R \in B^\pi$.
\end{lemma}
\smallskip
{\it Proof}. Let $\fqk$ be a left coisotropic subgroup. 
If $A=\sum_i e_i\otimes A_i\in \indr$, then we have 
that $(\pi\otimes\id)\Delta A_i = \sum_j a_{ij} \otimes A_j$. The proof relies
on the following identity
\begin{equation}
\label{eq_sesq1}
\sum_j \sum_{(A_j)} S^{-1}({A_j}_{(1)})a_{ij} \otimes {A_{j}}_{(2)} = \pi(1)\otimes A_i\;.
\end{equation}
To prove (\ref{eq_sesq1}) let's define $\epsilon_L:\fqg\rightarrow\fqk$, with 
$\epsilon_L(a) = \sum_{(a)} S^{-1}(a_{(2)})\pi(a_{(1)}) = \epsilon(a) \pi(1)$.
Then $(\epsilon_L\otimes\id)\Delta A_i=\pi(1)\otimes A_i$. Defining  
$\mu_L:\fqk\otimes\fqg\rightarrow\fqk$, $\mu_L(a\otimes b)= ba$, we have
$$
\begin{array}{rcl}
(\epsilon_L\otimes\id) \Delta A_i &=& (\mu_L\otimes\id) (\pi\otimes S^{-1}\otimes\id)
(\id\otimes\Delta) \Delta A_i  \cr
& = & (\mu_L\otimes\id)(\id\otimes S^{-1}\otimes\id)(\id\otimes\Delta)
\sum_j a_{ij}\otimes A_j \cr
& = & \sum_j \sum_{(A_j)}S^{-1}({A_j}_{(1)})a_{ij} \otimes {A_j}_{(2)}\;,
\end{array}
$$
from which we get the result. Finally,
$$
\begin{array}{rcl}
(\pi\otimes\id)\Delta\langle A,B \rangle_L &=& \sum_i(\pi\otimes\id) \Delta(A_i^*B_i) \cr
& = & \sum_i \Delta A_i^*(\pi\otimes\id)\Delta B_i \cr
& = & \sum_{ij} \sum_{(A_i)}{A_i}_{(1)}^* a_{ij} \otimes {A_i}_{(2)}^* B_j  \cr
& = & \sum_{ij} \sum_{(A_i)} \tau_K(S^{-1}({A_i}_{(1)})a_{ji}) \otimes {A_i}_{(2)}^* B_j  \cr
& = & \sum_j \pi(1)\otimes A_j^* B_j = \pi(1) \otimes \langle A,B\rangle_L\,.
\end{array}
$$

The analogous proof holds when $\fqk$ is a right coisotropic subgroup. 
If $A=\sum_i A_i\otimes e_i\in\indl$ then
$(\id\otimes\pi)\Delta A_i=\sum_j A_j\otimes b_{ij}$. By defining 
$\epsilon_R(a) = \sum_{(a)} \pi(a_{(2)}) S^{-1}(a_{(1)}) = \epsilon(a)\pi(1)$, and by applying 
$(\id\otimes\epsilon_R)\Delta$ on $A_i$ we arrive to 
\begin{equation}
\label{eq_sesq2}
\sum_j\sum_{(A_j)} {A_j}_{(1)}\otimes b_{ji} S^{-1}({A_j}_{(2)}) = A_i\otimes\pi(1)\;,
\end{equation}
from which
$$
\begin{array}{rcl}
(\id\otimes\pi)\Delta\langle A,B\rangle_R &=& \sum_{ij} (B_j\otimes b_{ji}) \Delta A_i^* \cr
&=& \sum_{ij} \sum_{(A_i^*)} B_j{A^*_{i}}_{(1)} \otimes\tau_K(b_{ij}S^{-1}({A_i}_{(2)})) \cr
&=& \sum_j B_j A^*_j\otimes\pi(1) = \langle A,B\rangle_R \otimes\pi(1) \;. \fidi
\end{array}
$$ 

\medskip
Let's now fix a left coisotropic subgroup $\fqk$ and let $B_\pi$ be the 
corresponding left homogeneous space.
Let $h:B_\pi\rightarrow\CC$ be a quasi invariant functional according to 
(\ref{def_phi}). Then, if $A,B\in\indr$, 

\begin{equation}
\label{scalare}
\langle A,B\rangle = h(\langle A,B\rangle_L)
\end{equation}
defines a sesquilinear form on $\indr$.
If $h$ is only quasi invariant and not invariant ({\it i.e.}$\phi[X]\neq\epsilon(X)$) the
induced representation is not unitary. We can nevertheless define a unitary representation.
Let $\rin$ be the representation dual to the induced corepresentation on $\indr$.

\medskip
\begin{proposition}
\label{pro_un_ind}
The application $\urin(X)$ defined in the basis $\{e_i\}$ of $V$ by
$$
\urin(X) A = \sum_{(X)} \sum_i e_i \otimes X_{(1)}.A_i\phi[X_{(2)}]
$$
is a unitary left representation of $\uq$ on $\indr$ with respect to the 
sesquilinear form {\rm(\ref{scalare})}, i.e. 
$$
\langle A,\urin(X) B \rangle = \langle \urin(X^*) A, B\rangle \;.
$$

Let $h_1$ and $h_2$ be equivalent quasi-invariant functionals on 
$B_{\pi}$ and let $\urinn 1$ and $\urinn 2$ be the corresponding induced
unitary representations. Then there is a unitary equivalence between $\urinn 1$
and $\urinn 2$.
\end{proposition}

\smallskip
{\it Proof}. From Lemma \ref{lem_B_mod} we see that $\urin$ maps $\indr$ into
$\indr$ and from (\ref{def_phi}) we verify that it is a representation of $\uq$. 

We prove its unitarity. Using basic properties of coproduct and antipode, we can write 
$ \langle A,\urin(X) B\rangle =\sum_i\sum_{(X)} 
h( X_{(2)}.(S^{-1}(X_{(1)}).A_i^* B_i) \phi[X_{(3)}]) $ and applying Lemma
\ref{lem_q_inv_m} we obtain 
$$
\langle A,\urin(X) B\rangle = \sum_i\sum_{(X)} 
h(\phi[X_{(2)}^*]^*S^{-1}(X_{(1)}).A_i^*B_i) \;.
$$ 
Since
$$
\langle \urin(X^*)A,B\rangle = \sum_i\sum_{(X)} 
h(\phi[X_{(2)}^*]^*(X_{(1)}^*.A_i)^* B_i)\;
$$
the result follows from the property $(X^*.a)^*=S^{-1}(X).a^*$. 

Concerning the second statement, let $h_1(a)=h_2(\xi^*a\xi)$ and 
let $F:\indl\rightarrow\indl$ be defined by 
$F(v\tens a)=v\tens a\xi$.
If $A \in \indl$, then
$$
\begin{array}{rcl}
F\circ \urinn 1 (X)A &=& \sum_i e_i\tens X_{(1)}.A_i \phi_1[X_{(2)}]\xi \\
                      &=& \sum_i e_i\tens X_{(1)}.(A_i\xi)\phi_2[X_{(2)}] 
= \urinn 2(X)\circ F(A)
\end{array}
$$
where we used (\ref{w_equiv}). By direct computation this equivalence
is proved to be unitary.\fidi

\medskip
The analogous property is true for a right coisotropic subgroup.
Let $h:B^\pi\rightarrow\CC$ a quasi invariant 
functional on the right homogeneous space. If $A,B\in\indl$ then
\begin{equation}
\label{scalare2}
\langle A,B \rangle = h(\langle A, B \rangle_R)
\end{equation}
defines a sesquilinear form on $\indl$. We have:

\medskip
\begin{proposition}
\label{pro_un_ind2}
The application $\ulin(X)$ defined by
$$
\ulin(X) A = \sum_{(X)} \sum_i \psi[X_{(1)}]A_i.X_{(2)}\otimes e_i
$$
is a unitary right representation of $\uq$ on $\indl$ with respect to the 
sesquilinear form {\rm(\ref{scalare2})}, i.e. 
$$
\langle A,\ulin(X) B \rangle = \langle \ulin(X^*) A, B\rangle \;.
$$

Let $h_1$ and $h_2$ be equivalent quasi-invariant functionals on $B^{\pi}$. 
and let $\ulinn 1$ and $\ulinn 2$ be the corresponding induced
unitary representations. Then there is a unitary equivalence between $\ulinn 1$
and $\ulinn 2$.
\end{proposition}

\bigskip
\bigskip

\sect{The Quantum Galilei Group with Imaginary Parameter}
The representation theory of the one dimensional quantum Galilei group with
real parameter was studied in \cite{F1}. In that case physical
interpretation required a non standard involution ({\it i.e.} $\tau^2\neq\id$).
The corresponding results on the unitarity of induced representations cannot be
treated in the framework just described. On the contrary the case with 
imaginary parameter and standard involution is a good example.

The one dimensional Quantum Galilei Group can be introduced by giving its
universal enveloping algebra $\uq(G(1))$, which is generated by $\{K,K^{-1},B,M\}$
and by the relations
$$
[K,T]=0\,,~~~~~~K\,B\,K^{-1}=B-iwM\,,~~~~~~[B,T]=i\fraz{K-K^{-1}}{2w}\,,
$$ 
where $M$ is central and $KK^{-1}=1=K^{-1}K$. The deformation parameter is $iw$
with real $w$. The coalgebra is given by

\begin{eqnarray*}
\Delta M &=&M\otimes K+K^{-1}\otimes M\,,~~~~~~~~~~
\Delta K=K\otimes K\,,\\
\Delta T &=&T\otimes 1 +1\otimes T \,,~~~~~
\Delta B=B\otimes K+K^{-1}\otimes B \,
\end{eqnarray*}
and
$$
S(M)=-M\,,~~~~~~~S(K)=K^{-1}\,,~~~~~~~S(T)=-T\,,~~~~~~~S(B)=-B+iwM\,;\\
$$
$$
\epsilon(M) =0\,, ~~~~~~~ \epsilon(K)=1\,, ~~~~~~~~\epsilon(T)=0\,, ~~~~~~\epsilon(B)=0\;.
$$
The Hopf algebra of polynomial functions $\fq(G(1))$ is generated by
$\{\mu,x,v,t\}$ and by the relations
$$
[\mu,x]=2iw\mu\,,~~~~~~~[\mu,v]=-iwv^2\,,~~~~~~~[x,v]=-2iwv\,,
$$
$t$ being a central element. The structure is completed by 

\begin{eqnarray*}
\Delta\mu &=& \mu\otimes 1+1\otimes\mu+v\otimes x +\half v^2\otimes t\,,~~~~
\Delta t=t\otimes 1+1\otimes t\,,\cr
\Delta x&=&x\otimes 1+1\otimes x+v\otimes t\,,~~~~~~~~
\Delta v =v\otimes 1+1\otimes v\,;
\end{eqnarray*}
$$
S(\mu)=-\mu+vx-\half v^2t\,,~~~~~~ S(x)=-x+tv\,,~~~~~S(t)=-t\,,
~~~~~S(v)=-v\,;
$$ 
$$
\epsilon(\mu)=0\,, ~~~~~~ \epsilon(x) =0\,,~~~~~~\epsilon(t)=0\,,~~~~~~
\epsilon(v)=0\,.
$$

The involution is given by defining all the generators of $\uq(G(1))$ and $\fq(G(1))$ 
to be real, with the exception of $\mu$ for which $\mu^*=\mu-iw v$.
Finally, the duality between $\uq(G(1))$ and $\fq(G(1))$ reads
$$
I^{\alpha'} K^{\ell} T^{\gamma'} N^{\delta'}(\mu^\alpha x^\beta t^\gamma v^\delta\,) =
i^{\alpha+\beta+\gamma+\delta}\ \alpha!\,\gamma!\,\delta!\,
(-iw\,\ell)^\beta\,\,\delta_{\alpha,\alpha'}\,
\delta_{\gamma,\gamma'}\,\delta_{\delta,\delta'}
$$
where $I=K^{-1} M,\, N=K B$ and $\ell\in {\bf Z}$ while the other
indices are in ${\bf N}$. The involution satisfies the usual property
$\tau^2=\id$ and $X^*(a)=\overline{X(\tau(a))}$.
 
It is easy to verify that the three primitive and real generators 
$\{\muh,\xh,\th\,\}$ with relations $[\muh,\xh]=2iw\muh$ and $[\th,\cdot]=0$
and the projection $\pi(\mu)=\muh,\pi(t)=\th,\pi(x)=\xh,\pi(v)=0$ define
a real quantum subgroup $\ffq J$.
Finally, $\omega_{m,u}=\exp[- i(m\muh+u\th\,)]$ with $m,u\in\R$ defines a unitary 
one-dimensional corepresentation of $\ffq J$.

The space $B_\pi$ is formed by the polynomials in $v$ and it must be extended
in order to support the induced representation. 

The construction is the same described in \cite{F1} for the real parameter case
to which we refer for details.
The algebra $\hoirr$ of quantum square 
integrable functions on the homogeneous space is described in terms of two commuting
generators $v_0$ and $v_1$ with the relation
$$
wmv_0v_1 = v_1-v_0 \;,
$$
and $v_0^*=v_0$, $v_1^*=v_1$. An useful basis is given by 
$\{\chi_m^\ell\}_{\ell\in\ZZ}$, where $\chi_m=1+wmv_1$ and 
$\chi_m^{-1}=1-wmv_0$. We can define the following structure
of $\uq$-module $*$-algebra on $\hoirr$:
$$
B.\chi_m^\ell=iwm\,\ell\,\chi_m^{\ell+1}\;,
$$
and $K.f=f$, $T.f=M.f=0$ for each $f\in\hoirr$. This algebra can be realized
in the dual of $\uq$. Indeed let $\fmu = e^{-i(m\mu+ut)}$, then 
$v_0=v$ and $v_{1}=\fmu^{-1} v \fmu$. The induced representation space is then 
$\hirr= \fmu\hoirr$.

Let's define the functional $\nu_w:\hoirr\rightarrow{\bf C}$,
$$
\nu_w(1)= 1\ , \quad\quad
\nu_w(v_0^n)=\fraz 1{(wm)^n}\ , \quad\quad \nu_w(v_{1}^n)=\fraz
1{(-wm)^n}\;,
$$
i.e. $\nu_w(\chi_m^\ell)=\delta_{\ell,0}$; let the corresponding sesquilinear 
form on $\hirr$ be
$\langle a,b \rangle=\nu_w(a^*b)$ for $a,b\in\hirr$. Then, since
$\fmu^*\fmu=\chi_m$, we have that
\begin{equation}
\label{p_scalar}
\langle \fmu \chi_m^\ell,\fmu\chi_m^n\rangle =\nu_w(\chi_m^{\ell+n+1})=
\delta_{\ell+n+1,0}  \;.
\end{equation}

According to \cite{Ota} we can consider on $\hoirr$ the scalar product
\begin{equation}
\label{truescalar}
(\fmu \chi_m^\ell, \fmu \chi_m^n) =\delta_{\ell,n}.
\end{equation}
If we define $j:\hirr\rightarrow \hirr$ by
$j(\fmu \chi_m^\ell)= \fmu \chi_m^{-\ell-1}$, we see that
\begin{equation}
\label{jform}
\pair {\fmu \chi_m^\ell}{\fmu \chi_m^n}=(j(\fmu \chi_m^\ell), \fmu \chi_m^n)
\end{equation}
and $\langle,\rangle$ is a Minkowski form on the $j$-space $\hoirr$.

The quasi-invariance properties of $\nu_w$ are collected in the following 
proposition. 

\smallskip
\begin{proposition}
\label{weight}
Let $\phi:\uq\rightarrow\hoirr$ be defined by
$$
\phi[X] = \epsilon(X) + \sum_{n=1}  
c_n\, X(v^n) \,\chi_m^n \,,~~~~~~~~ 
c_n=\left(\frac{wm}{2}\right)^n\ \frac{(2n-1)!!}{n!}\;.
$$
Then $\nu_w$ is quasi-invariant with weight $\phi$ but not essentially invariant.
\end{proposition}

\smallskip
{\it Proof}. From the definition of $\phi$ it follows that the cocycle
property (\ref{def_phi}) is verified for $X=\{M,T,K\}$ and arbitrary $Y\in\uq$. 
>From a direct computation 
we see that (\ref{def_phi}) for $X=B$ is a consequence of the following
recurrence relation
$$
nc_n= wm \left(n-\frac{1}{2}\right) c_{n-1}\;.
$$
The result follows from the properties of
the coproduct of generators. The property (\ref{def_h}) of quasi invariance
for $h$ is trivial for $X=\{M,T,K\}$ and can be directly checked on $\chi_m^\ell$
for $X=B$. \fidi

\medskip
Define ${\widetilde\rho_{m,u}}(X)a = \sum_{(X)}X_{(1)}.a\,\phi[X_{(2)}]$, with
$a\in\hirr$. By direct computation and using Proposition \ref{pro_un_ind} we obtain
the final result.

\smallskip
\begin{proposition}
\label{gal_rap_un}
The representation ${\widetilde\rho_{m,u}}$ of $\uq(G(1))$, given by
$$
\begin{array}{ll}
{}&\rrmu(K^{\pm 1})\,\fmu\, \chi_m^\ell= \fmu\, \chi_m^{\ell\pm 1}\,,
\quad\quad
\rrmu(B)\,\fmu\, \chi_m^\ell = iwm\,(\ell+\frac{1}{2})\,\fmu\, \chi_m^\ell\,,\cr
{}&\rrmu(T)\,\fmu\, \chi_m^\ell = \fmu\, \chi_m^\ell\, (\,-\fraz
1{2w^2m}\,(2-\chi_m-\chi_m^{-1})+u)\,; 
\end{array}
$$
is unitary with respect to the form {\rm(\ref{p_scalar})}.
\end{proposition}

According to the definitions of \cite{Ota} $\rrmu$ is a $j$-representation
of the quantum Galilei group.

\end{document}